\documentclass[journal]{IEEEtran}
\usepackage{amssymb,amsfonts}
\usepackage{algorithmic}
\usepackage{graphicx}
\usepackage{textcomp}
\def\BibTeX{{\rm B\kern-.05em{\sc i\kern-.025em b}\kern-.08em
    T\kern-.1667em\lower.7ex\hbox{E}\kern-.125emX}}
\usepackage{enumitem}
\usepackage{float}
\usepackage[tbtags]{amsmath}
\usepackage[tbtags]{mathtools}
\usepackage{accents}
\usepackage{multicol}
\usepackage[ruled,vlined]{algorithm2e}
\usepackage{cite}
\usepackage{pgfplots}
\usepackage{array,multirow}
\usepackage{hhline}
\usepackage{arydshln}

\usepackage{color,soul}

\usepackage[hidelinks=true]{hyperref}
\usepackage{cite}

\usepackage{xcolor}
\hypersetup{
    colorlinks,
    linkcolor={blue},
    citecolor={blue},
    urlcolor={blue}
}

\newtheorem{assumption}{Assumption}

\newtheorem{remark}{Remark}

\begin{document}
\title{Covariance-Based Multiple-Impulse Rendezvous Design}
\author{Amir Shakouri,
\thanks{Amir Shakouri is with the Department of Aerospace Engineering, Sharif University of Technology, Tehran, Iran, e-mail: \href{mailto:a_shakouri@ae.sharif.edu}{a$\_$shakouri@ae.sharif.edu}.}
        Maryam Kiani, and Seid H. Pourtakdoust
\thanks{Maryam Kiani (corresponding author) and Seid H. Pourtakdoust are with the Center for Research and Development in Space Science and Technology, Sharif University of Technology, Tehran, Iran, e-mails: \href{mailto:kiani@sharif.edu}{kiani@sharif.edu}, \href{mailto:pourtak@sharif.edu}{pourtak@sharif.edu}.}}

\maketitle

\begin{abstract}
A novel trajectory design methodology is proposed in the current work to minimize the state uncertainty in the crucial mission of spacecraft rendezvous. The trajectory is shaped under constraints utilizing a multiple-impulse approach. State uncertainty is characterized in terms of covariance, and the impulse time as the only affective parameter in uncertainty propagation is selected to minimize the trace of the covariance matrix. Further, the impulse location is also adopted as the other design parameter to satisfy various translational constraints of the space mission. Efficiency and viability of the proposed idea have been investigated through some scenarios that include constraints on final time, control effort, and maximum thruster limit addition to considering safe corridors. The obtained results show that proper selection of the impulse time and impulse position fulfils a successful feasible rendezvous mission with minimum uncertainty. 
\end{abstract}


\section{Introduction}
\label{sec:I}

\IEEEPARstart{S}TATE estimation of a dynamic system in the Bayesian framework needs a dynamic system and a measurement model, both of which could be usually contaminated with uncertainties emanating from mis-modeling and/or measurement errors. The uncertainties are usually represented as white Gaussian noises whose covariances are a measure of uncertainty showing how the data are distributed with respect to their expected values \cite{in1}. In this sense, many researchers have exploited the trace of the covariance matrix as a measure of total system uncertainty and as a criterion to evaluate dynamic system performance \cite{in2}. Other methods for uncertainty quantization have also been discussed in the literature such as using the determinant of the covariance matrix \cite{in2p}. Kalman filter is an optimal minimum variance state estimation algorithm for a linear dynamic system with a linear measurement model \cite{in3}. In a linear time-invariant (LTI) dynamic system with a linear measurement model, error covariance matrix propagates in time independent of states knowledge. This property creates an advantage to control the uncertainty level at some key positions via proper selection of pertinent time instants.  For nonlinear orbital maneuvers, the reader can refer to \cite{in3p,in3pp} where the uncertainty propagation in Lambert's problem is discussed in details. 

There are various approaches to describe the spacecraft dynamics in a rendezvous mission due to different assumptions one can consider for perturbations and/or orbital characteristics \cite{in4}. The Clohessy-Wiltshire (CW) model \cite{in5} develops a LTI system of equations for  rendezvous between close spacecraft in near circular orbits. There are also a wide variety of measurements for the absolute and relative navigation of spacecraft in a rendezvous mission including line-of-sight \cite{in6}, angles-only data \cite{in7}, vision-based navigation \cite{in8}, global positioning system (GPS) \cite{in9}, and relative GPS \cite{in10}. 

The covariance propagation idea as a driver for spacecraft trajectory design has previously been investigated by Zimmer et al. \cite{in11} but for low-thrust missions. Li et al. \cite{in12} have implemented a multi-objective unconstrained optimization for an impulsive rendezvous mission and considered the effect of uncertainties as a part of the cost function. Rendezvous missions mostly come across several constraints such as maximum control effort, maximum thruster limit and maximum flight time \cite{in13}. Moreover, safe corridors for the approach ellipsoid (AE) or analogously the keep-out ellipsoid (KOE) are important factors in mission safety analysis \cite{in14}.

In this paper a covariance-based approach is proposed to design a multiple-impulse spacecraft rendezvous. For a similar problem, Li et al. \cite{in12} tried to minimize unconstrained cost functions consisting of uncertainty and energy requirements while the current study is aimed at solutions that minimize the uncertainty subject to various constraints. Impulse time instants and their corresponding locations are considered as design parameters in the proposed methodology that is based on the CW equations. In addition, a linear measurement model is taken for the relative position and velocity via GPS data. The proposed method is implemented to design a two- and a three-impulse rendezvous problem for two case studies with different control effort and thruster limit constraints. The maximum flight time and the safe corridors such as the AE and the KOE are among the other mission new constraints considered in these scenarios. In this sense, the key contributions of the present work that differentiates it against the existing researches include: 

\begin{enumerate}[label=\emph{(\roman*)}]
\item Covariance trace is considered as the only driver for a completely impulsive rendezvous trajectory design;
\item The impulse time instants are used for uncertainty minimization, while the impulse locations are considered to address the mission constraints;
\item Maximum control effort, maximum impulse level, maximum flight time as well as the safe corridors are simultaneously taken as the mission constraints.
\end{enumerate}

The rest of the paper is organized as follows: Section \ref{sec:II} presents the preliminaries in which the multiple-impulse rendezvous has been formulated in a transition matrix form and the impulses are derived as functions of impulse time and location. Moreover, propagation of the uncertainty covariance is discussed in Section \ref{sec:II}. Section \ref{sec:III} explains the trajectory design methodology to select the impulse time and location appropriately. Section \ref{sec:IV} provides numerical simulation and analysis for two and three-impulse rendezvous. Finally, Section \ref{sec:V} concludes the paper and offers some recommendations for future work. 

\section{Preliminaries}
\label{sec:II}

\subsection{Notations}
\label{sub:II-A}

This subsection briefly introduces the notations that are used throughout the paper. 

Let $\mathbb{M}^{m,n}$ denote the space of $m \times n$ real (or complex) matrices and $\mathbb{M}^n$ its square analog. Let $\mathbb{S}^n$ denote the space of $n\times n$ symmetric matrices. In addition, $\mathbb{R}^n$ represents the space of $n$-dimensional real vectors and $\mathbb{I}_n$ indicates the $n$-dimensional identity matrix. The symbol $\|\cdot\|$ denotes the Euclidean norm of a vector. For matrix $M\in\mathbb{M}^{m,n}$, we note its transpose by $M^T$, its inverse (if exists) by $M^{-1}$, its trace by $\text{tr}(M)$. For $\rho\in\mathbb{R}$ the symbol $[\rho]_{m\times n}$ is used to denote a $m\times n$ matrix with all entries equal to $\rho$. 

The operator $\text{Diag}(\cdot):\mathbb{R}^n\mapsto\mathbb{S}^n$ maps a vector into a diagonal matrix with the diagonal elements equal to the vector entries. The expected operator is shown by $\text{E}(\cdot)$ and a normally distributed random vector, $\boldsymbol r\in\mathbb{R}^n$, is referred to by 
$$\boldsymbol r\sim\mathcal{N}\left(\text{E}(\boldsymbol r),\text{E}\left\{\left[\boldsymbol r-\text{E}(\boldsymbol r)\right]\left[\boldsymbol r-\text{E}(\boldsymbol r)\right]^T\right\}\right)$$

The TS-centered RSW coordinate system is defined as follows: The $x$-axis (R) is along with the radius of the target spacecraft (TS) orbit, the $z$-axis (W) extends along the angular momentum vector of the TS orbit, and the $y$-axis (S) completes the right-handed system. 
\subsection{Multiple-Impulse Rendezvous Formulation}
\label{sub:II-B}

The relative motion of a chaser spacecraft (CS) with respect to a TS can be formulated by a set of linear differential equations under the following assumption: 
\begin{assumption}
\label{ass1}
Let $a_t\in\mathbb{R}$ denote the semi-major axis of the TS's orbit and $\boldsymbol r\in\mathbb{R}^3$ denote the relative position of the CS with respect to the TS in an arbitrary TS-centered coordinate system. The following assumptions are made: 

\begin{enumerate}[label=\emph{(\roman*)}]
\item The two-body gravitational force is governing and any perturbation is ignored. 
\item The TS is in a circular orbit.
\item $\|\boldsymbol r\|/a_t\ll1$.
\end{enumerate}
\end{assumption}

Let $x,y,z\in\mathbb{R}$ denote the elements of the position vector of the CS in a TS-centered RSW coordinate system, $\boldsymbol r=[x\quad y\quad z]^T$, and $n=\sqrt{\mu/a_t^3}$. Under Assumption \ref{ass1}, the following system is obtainable \cite{in5}: 
\begin{equation}
\label{eq:1}
\ddot{x}=3n^2x+2n\dot{y}
\end{equation}
\begin{equation}
\label{eq:2}
\ddot{y}=-2n\dot{x}
\end{equation}
\begin{equation}
\label{eq:3}
\ddot{z}=-n^2z
\end{equation} 

The state of the CS at $t=t_i\in[0,\infty)$ is defined by $\boldsymbol r_i\triangleq[x_i\quad y_i\quad z_i]^T\in\mathbb{R}^3$ and $\boldsymbol v_i\triangleq[\dot{x}_i\quad \dot{y}_i\quad \dot{z}_i]^T\in\mathbb{R}^3$. The velocity impulse vector at $t=t_i$ is defined by $\Delta\boldsymbol v_i\triangleq[\Delta\dot{x}_i\quad\Delta\dot{y}_i\quad\Delta\dot{z}_i]^T$. Therefore, according to ‎Eqs. \eqref{eq:1}--\eqref{eq:3}, equations of motion of  CS in transferring from $\boldsymbol r_i$ to $\boldsymbol r_{i+1}$ by an impulse of $\Delta\boldsymbol v_i$, in a time interval of $Δt_{i(i+1)}\triangleq t_{i+1}-t_i$, are as below
\begin{equation}
\label{eq:4}
\left\{
\begin{matrix}
\boldsymbol r_{i+1} \\ \boldsymbol v_{i+1}
\end{matrix}
\right\}
=
\left[
\begin{matrix}
\Phi_{rr}\left(\Delta t_{i(i+1)}\right) & \Phi_{rv}\left(\Delta t_{i(i+1)}\right) \\
\Phi_{vr}\left(\Delta t_{i(i+1)}\right) & \Phi_{vv}\left(\Delta t_{i(i+1)}\right) \\
\end{matrix}
\right]
\left\{
\begin{matrix}
\boldsymbol r_{i} \\ \boldsymbol v_{i}+\Delta\boldsymbol v_{i}
\end{matrix}
\right\}
\end{equation}
where $\Phi_{rr}(\cdot),\Phi_{rv}(\cdot),\Phi_{vr}(\cdot),\Phi_{vv}(\cdot):[0,\infty)\mapsto\mathbb{M}^3$ such that
\begin{equation}
\label{eq:5}
\Phi_{rr}(t)\triangleq
\left[
\begin{matrix}
4-3\cos nt & 0 & 0 \\
6(\sin nt-nt) & 1 & 0 \\
0 & 0& \cos nt
\end{matrix}
\right]
\end{equation}
\begin{equation}
\label{eq:6}
\Phi_{rv}(t)\triangleq\frac{1}{n}
\left[
\begin{matrix}
\sin nt & 2(1-\cos nt) & 0 \\
-2(1-\cos nt) & 4\sin nt-3nt & 0 \\
0 & 0& \sin nt
\end{matrix}
\right]
\end{equation}
\begin{equation}
\label{eq:7}
\Phi_{vr}(t)\triangleq n
\left[
\begin{matrix}
3\sin nt & 0 & 0 \\
-6(1-\cos nt) & 0 & 0 \\
0 & 0& -\sin nt
\end{matrix}
\right]
\end{equation}
\begin{equation}
\label{eq:8}
\Phi_{vv}(t)\triangleq
\left[
\begin{matrix}
\cos nt & 2\sin nt & 0 \\
-2\sin nt & 4\cos nt-3 & 0 \\
0 & 0& \cos nt
\end{matrix}
\right]
\end{equation}

\begin{figure}[!h]
\centering\includegraphics[width=0.65\linewidth]{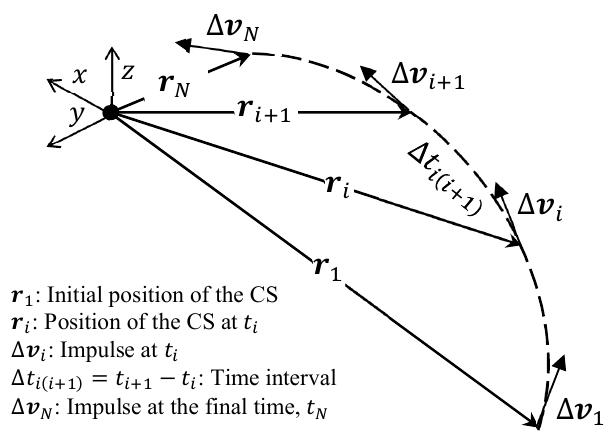}
\caption{Geometry of $N$-impulse rendezvous.}
\label{fig:1}
\end{figure}

As shown in Fig. \ref{fig:1}, through an $N$-impulse rendezvous maneuver, the CS should reach $\boldsymbol r_N$ and $\boldsymbol v_N$, starting from the initial state $\boldsymbol r_1$ and $\boldsymbol v_1$. A single impulse in the initial position, $\boldsymbol r_1$, $N-2$ impulses in the mid-positions, $\boldsymbol r_i$, $i=2,...,N-1$ and finally, a single impulse at the TS, $\boldsymbol r_N$ are needed to this aim. 

Assume that the velocity and position vectors of the CS at $t_1$ are previously known via navigation system. After specifying the mid-positions, $\boldsymbol r_i$ and the time intervals, $\Delta t_{i(i+1)}$, the impulse vectors can be obtained from Eqs. \eqref{eq:4} to \eqref{eq:8}: 
\begin{equation}
\label{eq:9}
\begin{split}
\Delta\boldsymbol v_i=-\Phi_{rv}^{-1}\left(\Delta t_{i(i+1)}\right)\Phi_{rr}\left(\Delta t_{i(i+1)}\right)\boldsymbol r_i \\
+\Phi_{rv}^{-1}\left(\Delta t_{i(i+1)}\right)\boldsymbol r_{i+1}-\boldsymbol v_i
\end{split}
\end{equation}
The velocity vector at $t_i$ can be obtained similarly as
\begin{equation}
\label{eq:10}
\begin{split}
\boldsymbol v_i=\left[\Phi_{vr}\left(\Delta t_{(i-1)i}\right)-\Phi_{vv}\left(\Delta t_{(i-1)i}\right)\Phi_{rv}^{-1}\left(\Delta t_{(i-1)i}\right)\right]\boldsymbol r_{i-1} \\
+\Phi_{vv}\left(\Delta t_{(i-1)i}\right)\Phi_{rv}^{-1}\left(\Delta t_{(i-1)i}\right)\boldsymbol r_i
\end{split}
\end{equation}
Substituting Eq. \eqref{eq:10} into Eq. \eqref{eq:9} yields the following expression for $i=1,...,N$:
\begin{equation}
\label{eq:11}
\begin{split}
\Delta\boldsymbol v_i=\Upsilon_1\left(\Delta t_{(i-1)i}\right)\boldsymbol r_{i-1}+\left[\Upsilon_2\left(\Delta t_{(i-1)i}\right) \right.\\
\left.+\Upsilon_3\left(\Delta t_{(i-1)i}\right)\right]\boldsymbol r_i+\Upsilon_4\left(\Delta t_{(i-1)i}\right)\boldsymbol r_{i+1}
\end{split}
\end{equation}
where $\Upsilon_1(\cdot),\Upsilon_2(\cdot),\Upsilon_3(\cdot),\Upsilon_4(\cdot):[0,\infty)\mapsto\mathbb{M}^3$ such that
\begin{equation}
\label{eq:12}
\Upsilon_1(t)\triangleq\Phi_{vv}(t)\Phi_{rv}^{-1}(t)\Phi_{rr}(t)-\Phi_{vr}(t)
\end{equation}
\begin{equation}
\label{eq:13}
\Upsilon_2(t)\triangleq-\Phi_{vv}(t)\Phi_{rv}^{-1}(t)
\end{equation}
\begin{equation}
\label{eq:14}
\Upsilon_3(t)\triangleq-\Phi_{rv}^{-1}(t)\Phi_{rr}(t)
\end{equation}
\begin{equation}
\label{eq:15}
\Upsilon_4(t)\triangleq\Phi_{rv}^{-1}(t)
\end{equation}
For the initial impulse, from Eq. \eqref{eq:9} we have:
\begin{equation}
\label{eq:16}
\Delta\boldsymbol v_1=\Upsilon_3(\Delta t_{12})\boldsymbol r_1+\Upsilon_4(\Delta t_{12})\boldsymbol r_2-\boldsymbol v_1
\end{equation}
where $\boldsymbol v_1$, and $\boldsymbol r_1$ are the CS initial position and velocity vectors. For the final impulse, 
\begin{equation}
\label{eq:17}
\Delta\boldsymbol v_N=\boldsymbol v_N^+-\boldsymbol v_{N}
\end{equation}
using Eq. \eqref{eq:10}, we have: 
\begin{equation}
\label{eq:18}
\begin{split}
\Delta\boldsymbol v_N=\boldsymbol v_N^++\Upsilon_1\left(\Delta t_{(N-1)N}\right)\boldsymbol r_{N-1} \\
+\Upsilon_2\left(\Delta t_{(N-1)N}\right)\boldsymbol r_N
\end{split}
\end{equation}
where $\boldsymbol v_N^+$ is the fixed velocity vector at $t_N$.

\subsection{Uncertainty Propagation}
\label{sub:II-C}

Uncertainties in the initial states of the CS and the inaccuracy of the system model are the main error sources that make the nominal trajectory deviate from the true trajectory. Further, the navigation errors due to GPS sensors degrade the performance of the measurement process. The covariance matrix, $P\in\mathbb{S}^6$, is a common measure to characterize the state uncertainty. The covariance matrix propagates in time under the influence of the process and  measurement  errors \cite{in15}. Initial covariance matrix, $P_1$, represents the uncertainty in initial state that could have emanated out of navigation error. The uncertainty in the system dynamics is modeled via a zero-mean Gaussian white noise, $\boldsymbol\omega\sim\mathcal{N}(0,Q)\in\mathbb{R}^6$ with $Q\triangleq⁡\text{E}(\boldsymbol\omega\boldsymbol\omega^T)\in\mathbb{S}^6$. The uncertainty in the measurements is also modeled by a zero-mean Gaussian white noise, $\boldsymbol\nu\sim\mathcal{N}(0,R)\in\mathbb{R}^6$, where $R\triangleq⁡\text{E}⁡(\boldsymbol\nu\boldsymbol\nu^T)\in\mathbb{S}^6$. It is assumed that the measurement noise covariance is time-invariant.

Let $F\in\mathbb{M}^6$ and consider the process model as follows: 
\begin{equation}
\label{eq:19}
\left\{\begin{matrix}
\dot{\boldsymbol r} \\
\dot{\boldsymbol v}
\end{matrix}\right\}=F
\left\{\begin{matrix}
\boldsymbol r \\
\boldsymbol v
\end{matrix}\right\}+\boldsymbol\omega
\end{equation}
where according to Eqs. \eqref{eq:1}--\eqref{eq:3}
\begin{equation}
\label{eq:20}
F\triangleq\left[\begin{matrix}
0 & 0 & 0 & 1 & 0 & 0  \\
0 & 0 & 0 & 0 & 1 & 0  \\
0 & 0 & 0 & 0 & 0 & 1  \\
3n^2 & 0 & 0 & 0 & 2n & 0  \\
0 & 0 & 0 & -2n & 0 & 0  \\
0 & 0 & -n^2 & 0 & 0 & 0  \\
\end{matrix}\right]
\end{equation}

Let $H\in\mathbb{M}^6$ and consider the measurement model to be
\begin{equation}
\label{eq:21}
\boldsymbol z=H\left\{\begin{matrix}
\boldsymbol r \\
\boldsymbol v
\end{matrix}\right\}+\boldsymbol \nu
\end{equation}
where $\boldsymbol z\in\mathbb{R}^6$ supposed to be the relative position and velocity of the CS obtained by a relative GPS \cite{in10}, i.e.,  $H\triangleq\mathbb{I}_6$. The state covariance matrix propagates in time for the above mentioned system as below \cite{in3},
\begin{equation}
\label{eq:22}
\begin{split}
P_{i+1}=\left[\Theta_{11}\left(\Delta t_{i(i+1)}\right)P_i+\Theta_{12}\left(\Delta t_{i(i+1)}\right)\right] \\
\times\left[\Theta_{21}\left(\Delta t_{i(i+1)}\right)P_i+\Theta_{22}\left(\Delta t_{i(i+1)}\right)\right]^{-1}
\end{split}
\end{equation}
where $\times$ denotes the conventional matrix multiplication, $\Theta_{11}(\cdot),\Theta_{12}(\cdot),\Theta_{21}(\cdot),\Theta_{22}(\cdot):[0,\infty)\mapsto\mathbb{M}^6$, and $\Lambda\in\mathbb{M}^{12}$ such that
\begin{equation}
\label{eq:23}
\left[\begin{matrix}
\Theta_{11}(t) & \Theta_{12}(t) \\
\Theta_{21}(t) & \Theta_{22}(t)
\end{matrix}\right]\triangleq\exp(\Lambda t)
\end{equation}
and
\begin{equation}
\label{eq:24}
\Lambda\triangleq\left[\begin{matrix}
F & Q \\
H^TRH & -F^T
\end{matrix}\right]
\end{equation}

Due to stationary properties of the above system, the steady state value of $P_i\equiv P_\infty$, at $t_i\rightarrow\infty$ can be obtained by the solution of the following algebraic Riccati equation (ARE) \cite{in3}: 
\begin{equation}
\label{eq:25}
FP_{\infty}+P_{\infty}F^T-P_{\infty}H^TR^{-1}HP_{\infty}+Q=0
\end{equation}

\begin{remark}
The uncertainty propagation formula discussed in this section can be used for any  linear time-variant (LTV) or LTI system represented in a state space form. Therefore, not only CW equations but the state-transition matrix for relative spacecraft motion developed in \cite{16} can be used as a basis for the trajectory design that is discussed in the next section. However, our analysis is based on the CW equations of motion. 
\end{remark}

\section{Trajectory Design Methodology}
\label{sec:III}

The main goal of the present study is to design a rendezvous trajectory that maximizes the probability of locating the spacecraft actual position at some selective key positions coincident with its nominal design positions. These selective positions are taken as the so called impulse locations in the current investigation. In this regard, the accuracy of the predicted trajectory is inserted into the rendezvous design procedure in terms of $\text{tr}(P_i)$. Besides, the total impulse and trajectory constraints including maximum accessible impulse value and the acceptable flight corridors are also considered as other design motivators. 

Therefore, in general the design parameters consist of the impulse time instants and their locations that create $4N-6$ adjustable parameters in an $N$-impulse rendezvous maneuver. The $N$ parameters of impulse time instants are determined mostly by analyzing the uncertainty behavior of the system while the impulse locations are selected to satisfy the problem constraints, where they may also cause a change in the control effort, that contain the remaining $3(N-2)$ design parameters.

\subsection{Impulse Time Instants}
\label{sub:III-A}

Since the process and the measurement models are linear, the covariance propagation can be accomplished offline. In this sense, the impulse time as our key parameter is chosen in accordance with  the time history of $\text{tr}(P_i)$. The desired covariance trace at $t_{i+1}$ can be chosen in the codomain of Eq. \eqref{eq:22}, from $\text{tr}(P_i)$ to $\text{tr}(P_\infty)$. 

The determined impulse time instants may not be compatible with the design constraints. In addition, the remaining design parameters such as the impulse positions and the problem constraints could also be dis-satisfied by any choice of impulse locations for some inconsistent impulse time.

A maximum possible final time may be determined according to the mission requirements. In this case, if the covariance has a converging behavior (that is regular in the presence of measurements and system observability), the final time, $t_N$, should be set at its maximum value corresponding to the minimum uncertainty at the destination. As the covariance is the key driver of the current design method, the impulse time instants should be preset in order to minimize the covariance trace until the problem constraints are satisfied. 

It is worth mentioning that the first impulse does not necessarily need to be exerted at $t=0$. In a system in which the covariance has a converging behavior, the first impulse can be applied with a delay in order to achieve lower uncertainty in the key location at $t_1$. This delay should be calculated considering the safe corridor. 

\subsection{Impulse Positions}
\label{sub:III-B}

The impulse control vectors are assumed to be applied with no error, thus their locations do not affect the covariance time history and the impulse values are functions of impulse locations. Therefore, the impulse locations are the remaining design variables needed to satisfy the mission constraints. Accordingly, the acceptable regions for impulse locations are determined based on the maximum allowable control effort, the upper thruster limit as well as the acceptable zones in space for the present study. 

For a maximum acceptable control effort, $J_{\text{max}}$, the following upper bound should be satisfied:

\begin{equation}
\label{eq:27}
\sum_{i=1}^N\|\Delta\boldsymbol v_i\|\leq J_{\text{max}}
\end{equation}
so, from Eqs. \eqref{eq:11}, \eqref{eq:16}, and \eqref{eq:18}, the impulse locations, $\boldsymbol r_i$, $i=2,...,N-1$, are selected in order to be inside the following hyper-surface: 

\begin{equation}
\label{eq:28}
\begin{split}
& \|\Upsilon_3\left(\Delta t_{12}\right)\boldsymbol r_1+\Upsilon_4\left(\Delta t_{12}\right)\boldsymbol r_2-\boldsymbol v_1\| \\
& +\sum_{i=2}^{N-1}\left\{\|\Upsilon_1\left(\Delta t_{(i-1)i}\right)\boldsymbol r_{i-1} \right.\\
& +\left[\Upsilon_2\left(\Delta t_{(i-1)i}\right)+\Upsilon_3\left(\Delta t_{(i-1)i}\right)\right]\boldsymbol r_i \\
& \left.+\Upsilon_4\left(\Delta t_{(i-1)i}\right)\boldsymbol r_{i+1}\|\right\} \\
& +\|\Upsilon_1\left(\Delta t_{(N-1)N}\right)\boldsymbol r_{N-1} \\
& +\Upsilon_2\left(\Delta t_{(N-1)N}\right)\boldsymbol r_N\|\leq J_{\text{max}}
\end{split}
\end{equation}

The maximum thruster limit $\Delta v_{\text{max}}$, poses another restriction to the problem as follows,
\begin{equation}
\label{eq:29}
\|\Delta\boldsymbol v_i\|\leq \Delta v_{\text{max}},\quad i=1,...,N
\end{equation}
hence, Eq. \eqref{eq:29} defines another bound for $\boldsymbol r_i$, $i=2,...,N-1$. Moreover, for safety reasons such as collision avoidance and cooperative operations, some space zones may be inaccessible. Suppose that the trajectory safe corridor is defined outside the KOE of $\boldsymbol r^TL_1\boldsymbol r=1$ and inside the AE of $\boldsymbol r^T L_2\boldsymbol r=1$, both of which are centered at the TS. Therefore: 
\begin{equation}
\label{eq:30}
\boldsymbol r^T(t)L_1\boldsymbol r(t)\geq1, \quad t=[t_1,t_N]
\end{equation}
\begin{equation}
\label{eq:31}
\boldsymbol r^T(t)L_2\boldsymbol r(t)\leq1, \quad t=[t_1,t_N]
\end{equation}
in which, $L_j=\text{Diag}\left\{\left[1/a_j^2\quad1/b_j^2\quad1/c_j^2\right]\right\}$ defines the semi-principal axes lengths of the ellipsoids. Therefore, Eqs. \eqref{eq:30} and \eqref{eq:31} using Eqs. \eqref{eq:4}, \eqref{eq:10}, \eqref{eq:11}, and \eqref{eq:16} now reduces to the following form: 
\begin{equation}
\label{eq:32}
\begin{split}
& \boldsymbol r_i^T\left[\Phi_{rr}(t)+\Phi_{rv}(t)\Upsilon_3\left(\Delta t_{i(i+1)}\right)\right]^TL_j\left[\Phi_{rr}(t)\right. \\
& \left.+\Phi_{rv}(t)\Upsilon_3\left(\Delta t_{i(i+1)}\right)\right]\boldsymbol r_i \\
& +2\boldsymbol r_{i+1}^T\Upsilon_4^T\left(\Delta t_{i(i+1)}\right)\Phi_{rv}^T(t)L_j\left[\Phi_{rr}(t)\right. \\
& \left.+\Phi_{rv}(t)\Upsilon_3\left(\Delta t_{i(i+1)}\right)\right]\boldsymbol r_i  \\
& +\boldsymbol r_{i+1}^T\Upsilon_4^T\left(\Delta t_{i(i+1)}\right)\Phi_{rv}^T(t)L_j\Phi_{rv}(t)\Upsilon_4\left(\Delta t_{i(i+1)}\right)\boldsymbol r_{i+1} \\
& \geq 1 
\end{split}
\end{equation}

Thus, from $t_i$ to $t_{i+1}$, those impulse locations that are applied in the region defined by Eq. \eqref{eq:31}, and lead the trajectory of the CS to satisfy Eq. \eqref{eq:32} for $j=1$ and $2$ are acceptable. Inequality \eqref{eq:32} guarantees that the CS does not deviate from the safe corridor after the $i$th impulse up to the next one. 

\begin{figure}[!h]
\centering\includegraphics[width=1\linewidth]{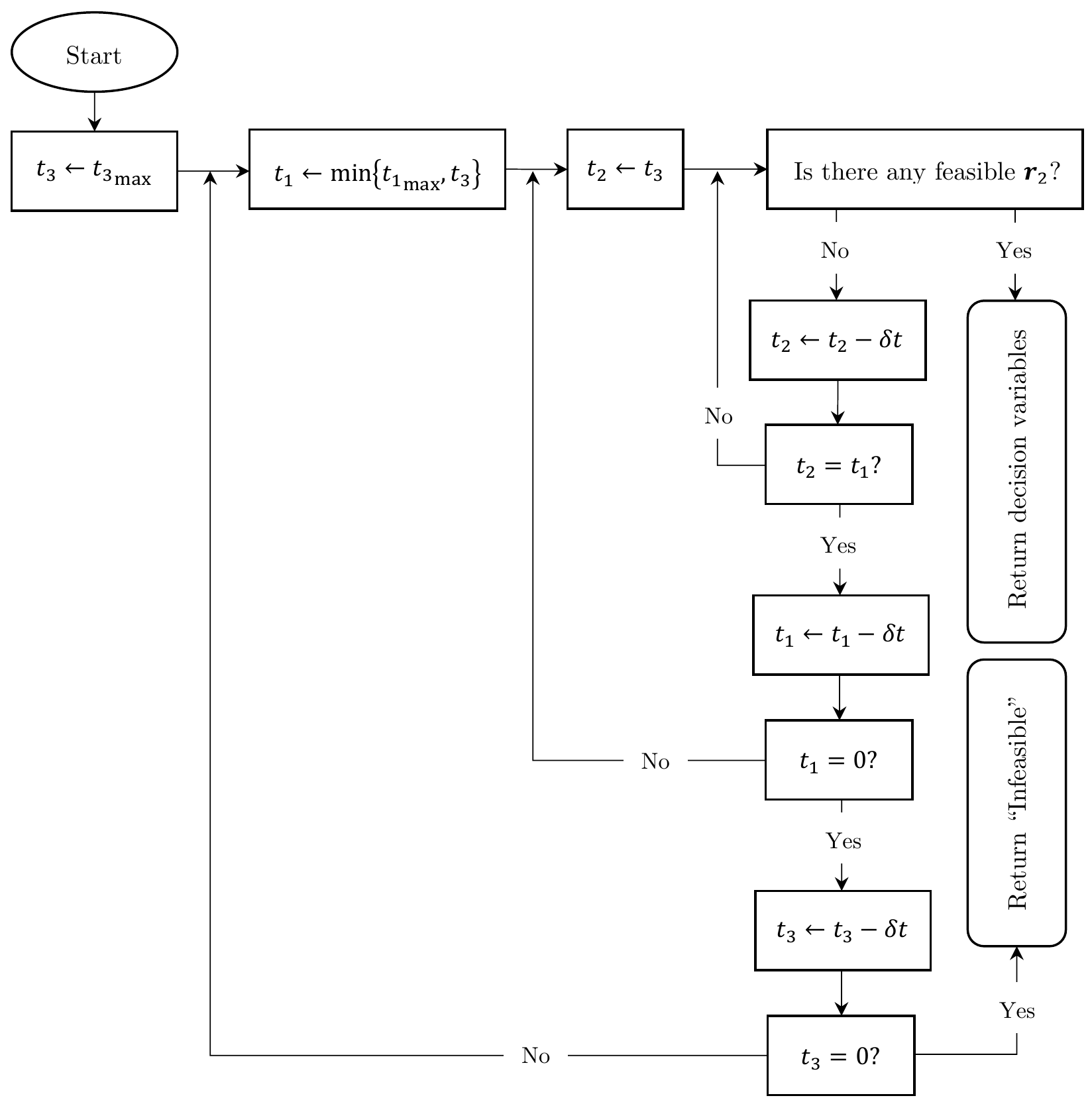}
\caption{A flow diagram representation of Algorithm \ref{al2}.}
\label{fig:2p}
\end{figure}

Algorithms \ref{al1} and \ref{al2} make use of the above-mentioned design procedure and search for solutions for cases of two- and three-impulse rendezvous problems. Step 2 in these algorithms utilizes Eqs. \eqref{eq:27}--\eqref{eq:32} to check the satisfaction of the constraints that are related to the impulse positions. The possibility of infeasible solutions is also considered in the algorithms. Fig. \ref{fig:2p} shows a flow diagram representation of Algorithm \ref{al2}. 

\begin{algorithm}[h]
 \KwData{Initial position and velocity vectors ($\boldsymbol r_1$ and $\boldsymbol v_1$), final position and velocity vectors ($\boldsymbol r_2$ and $\boldsymbol v_2$), table of constraints, and a time step ($\delta t>0$).}
 \KwResult{Time of impulses ($t_1^*$ and $t_2^*$).}
 \nl Select $t_1$ and $t_2$ (or $\Delta t_{12}$) at their maximum possible values, $t_1={t_1}_\text{max}$ and $t_2={t_2}_\text{max}$ (or $\Delta t_{12}={t_2}_\text{max}-{t_1}_\text{max}$). The value of ${t_2}_\text{max}$ is given as a constraint and the value of ${t_1}_\text{max}$ is determined at its maximum value such that the constraints are satisfied.
 
 \For{$t_2={t_2}_\text{max}:-\delta t:0$}{
 \For{$t_1=\min⁡{\{{t_1}_\text{max},t_2\}}:-\delta t:0$}{
  \nl If the constraints are satisfied break the loop, else continue.
  }
 }
 
 \If {$t_2=0$}{
 \Return{Infeasible}
 }
 \If{$t_2\neq0$}{
 \Return{$t_1^*=t_1$ and $t_2^*=t_2$}.}
 \caption{Covariance-based two-impulse rendezvous design algorithm.}
 \label{al1}
\end{algorithm}

\begin{algorithm}[h]
 \KwData{Initial position and velocity vectors ($\boldsymbol r_1$ and $\boldsymbol v_1$), final position and velocity vectors ($\boldsymbol r_3$ and $\boldsymbol v_3$), table of constraints, and a time step ($\delta t>0$).}
 \KwResult{Time of impulses ($t_1^*$, $t_2^*$, and $t_3^*$) and position of the second impulse ($\boldsymbol r_2^*$).}
 \nl Select $t_1$, $t_2$, and $t_3$ at their maximum possible values, $t_1={t_1}_\text{max}$ and $t_2=t_3={t_2}_\text{max}$. The value of ${t_2}_\text{max}$ is given as a constraint and the value of ${t_1}_\text{max}$ is determined at its maximum value such that the constraints are satisfied.
 
 \For{$t_3={t_3}_\text{max}:-\delta t:0$}{
  \For{$t_1=\min⁡{\{{t_1}_\text{max},t_3\}}:-\delta t:0$}{
   \For{$t_2=t_3:-\delta t:t_1$}{
    \nl If any pair of ${r_2}_x$ and ${r_2}_y$ can be found such that the constraints are satisfied break the loop, else continue. 
   }
  }
 }
 
 \If {$t_3=0$}{
 \Return{Infeasible}
 }
 \If{$t_3\neq0$}{
 \Return{$t_1^*=t_1$, $t_2^*=t_2$, $t_3^*=t_3$, and $\boldsymbol r_2^*=[{r_2}_x\quad {r_2}_y]^T$}.}
 \caption{Covariance-based three-impulse rendezvous design algorithm.}
 \label{al2}
\end{algorithm}

\section{Simulation Results}
\label{sec:IV}

Two case studies with different mission constraints are analyzed for a coplanar rendezvous with the following initial and final conditions expressed in the TS-centered RSW coordinate system:
$$\boldsymbol r_1=-[10^3\quad10^3\quad0]\text{ m},\quad\boldsymbol v_1=[10\quad-10\quad0]\text{ m/s}$$
$$\boldsymbol r_3=-[10^2\quad0\quad0]\text{ m},\quad\boldsymbol v_3=[0\quad0\quad0]\text{ m/s}$$

The TS is located in a circular-equatorial orbit with an altitude of $400\text{ km}$. The initial state covariance is assumed as
$$P_1=\text{Diag}\left\{\left[10^2[1]_{1\times3}\text{ m}^2\quad[1]_{1\times3}\text{ m}^2/\text{s}^2\right]^T\right\}$$
the covariance of the process noise is taken as $$Q=\text{Diag}\left\{\left[[0]_{1\times3}\text{ m}^2/\text{s}^2\quad10^{-12}[1]_{1\times3}\text{ m}^2/\text{s}^4\right]^T\right\}$$
and the covariance of the measurement noise is taken as
$$R=\text{Diag}\left\{\left[[1]_{1\times3}\text{ m}^2\quad10^{-2}[1]_{1\times3}\text{ m}^2/\text{s}^2\right]^T\right\}$$
The system constraints are supposed to be as shown in Table \ref{table:1} for both case studies.

\begin{table}[ht]
\caption{Constraints of the rendezvous case studies.}
\begin{center}
\begin{tabular}{|c|c|c|}
	\hline
	\multirow{2}{*}{Parameter}&\multicolumn{2}{c|}{Constraints} \\
	\cline{2-3}
	\multicolumn{1}{|c|}{} & Case Study 1 & Case Study 2 \\
	\hline
	Control Effort & $J_\text{max}=30\text{ m/s}$ & $J_\text{max}=20\text{ m/s}$ \\
	\hline
	Impulse Norm & $\Delta v_\text{max}=20\text{ m/s}$ & $\Delta v_\text{max}=13\text{ m/s}$ \\
	\hline
	Rendezvous Time & \multicolumn{2}{c|}{${t_N}_\text{max}=2\text{ h}$} \\
	\hline
	KOE & \multicolumn{2}{c|}{$a_1=50\text{ m}, \quad b_1=60\text{ m}, \quad c_1=70\text{ m}$} \\
	\hline
	AE & \multicolumn{2}{c|}{$a_2=10^4\text{ m}, \quad b_2=10^4\text{ m}, \quad c_2=10^4\text{ m}$} \\
	\hline
\end{tabular}
\end{center}
\label{table:1}
\end{table}

\begin{figure}[!h]
\centering\includegraphics[width=0.8\linewidth]{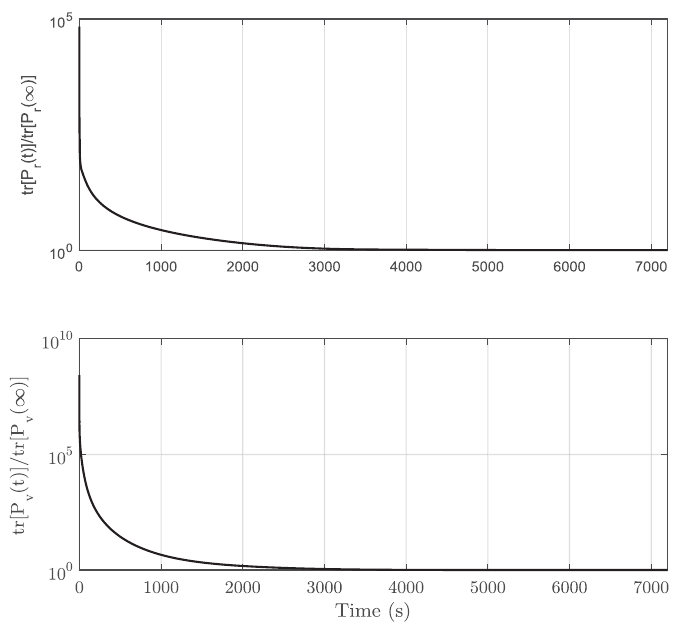}
\caption{Logarithmic diagram of normalized estimation performance of location (top) and velocity (down).}
\label{fig:2}
\end{figure}

The state covariance matrix is defined in a block form as follows: 
\begin{equation}
\label{eq:33}
P\triangleq\left[\begin{matrix}
P_r & P_{rv} \\
P_{rv}^T & P_{v}
\end{matrix}\right]
\end{equation}
The measurement units for $P_r$ and $P_v$ are $\text{m}^2$ and $\text{m}^2/\text{s}^2$, respectively. The normalized time history of the predicted uncertainty for position ($\text{tr}\left[P_r(t)\right]/\text{tr}\left[P_r(\infty)\right]$) and velocity ($\text{tr}\left[P_v(t)\right]/\text{tr}\left[P_v(\infty)\right]$) are shown in Fig. \ref{fig:2}. The stationary covariance values are obtained via the ARE equation, given in Eq. \eqref{eq:25}, and are equal to $\text{tr}\left[P_r(\infty)\right]=4.4169\times10^{-3}$ and $\text{tr}\left[P_v(\infty)\right]=1.1170\times10^{-8}$. Initially, a two-impulse transfer is studied for two case studies and then the same problem is carried out for a three-impulse scenario. The results have been summarized and discussed in Section \ref{sub:IV-C}.

According to Fig. \ref{fig:2} the uncertainty  extent can significantly change in the time domain (namely from $10^5$ and $10^{10}$ to $1$ for the position and velocity vectors, respectively). This converging behavior is more severe for poor initial estimations subject to more accurate observations. Therefore, it can be found that a covariance-based design can enhance the accuracy and success of a mission, while ignoring the covariance behavior would expose the system in risk. 

\subsection{Two-Impulse Rendezvous Design}
\label{sub:IV-A}

\begin{figure}[!h]
\centering\includegraphics[width=0.8\linewidth]{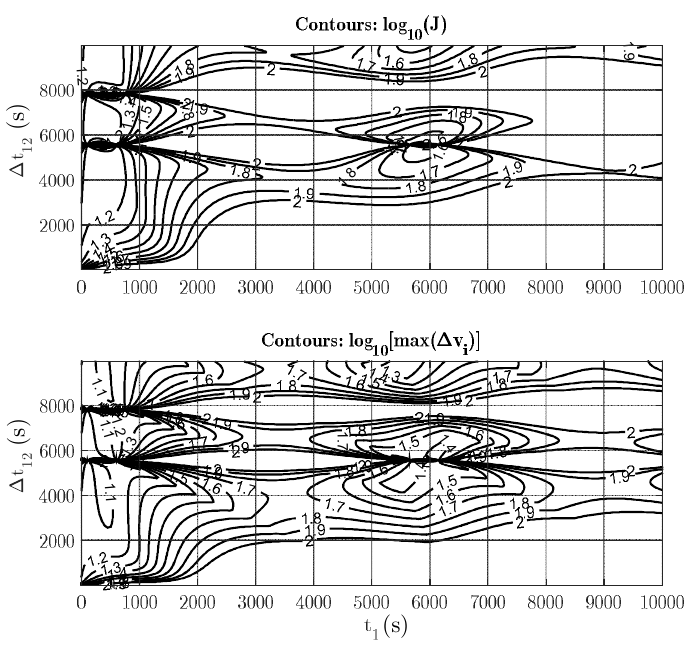}
\caption{Logarithmic contour plots for control effort (top) and maximum impulse (down) against $t_1$ and $\Delta t_{12}$.}
\label{fig:3}
\end{figure}

\begin{figure}[!h]
\centering\includegraphics[width=1\linewidth]{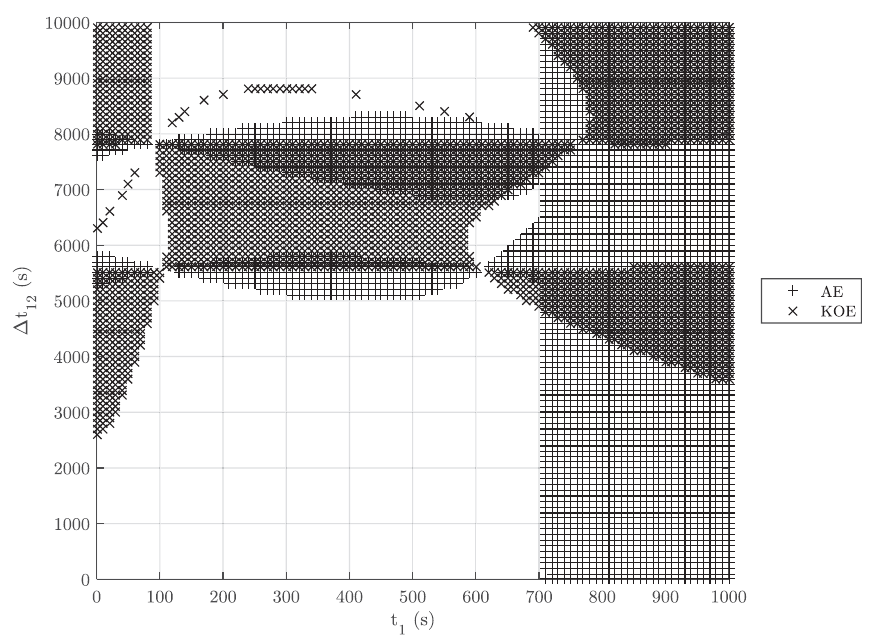}
\caption{Trajectories at which the AE and KOE are violated are shown by $+$ and $\times$ symbols, respectively.}
\label{fig:4}
\end{figure}

\begin{figure}[!h]
\centering\includegraphics[width=1\linewidth]{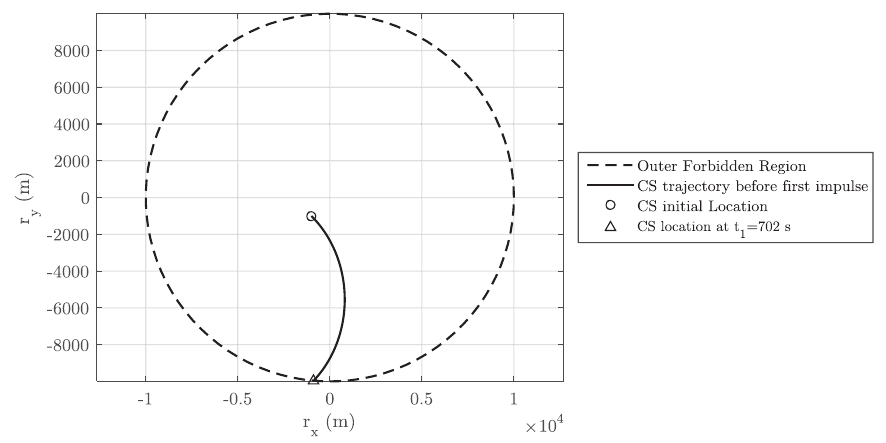}
\caption{CS trajectory from initial location to the intersection of the outer forbidden region.}
\label{fig:5}
\end{figure}

For a two-impulse rendezvous problem, the design parameters are obviously $t_1$ and $t_2$. According to Fig. \ref{fig:2} the minimum achievable uncertainty occurs for the maximum allowable rendezvous time, ${t_2}_\text{max}$. Fig. \ref{fig:3} shows the logarithmic contour plots of the control effort and the maximum impulse level against $t_1$ and $\Delta t_{12}$. Fig. \ref{fig:4} shows the time instants at which the associated trajectories violate the safety corridors of KOE and/or AE. At $t_1=702\text{ s}$ the location of the first impulse touches the AE (Fig. \ref{fig:5}) and at $t_1>702\text{ s}$ the trajectory violates the AE for every choice of $\Delta t_{12}$. 

\subsubsection*{Case Study 1}

\begin{figure}[!h]
\centering\includegraphics[width=1\linewidth]{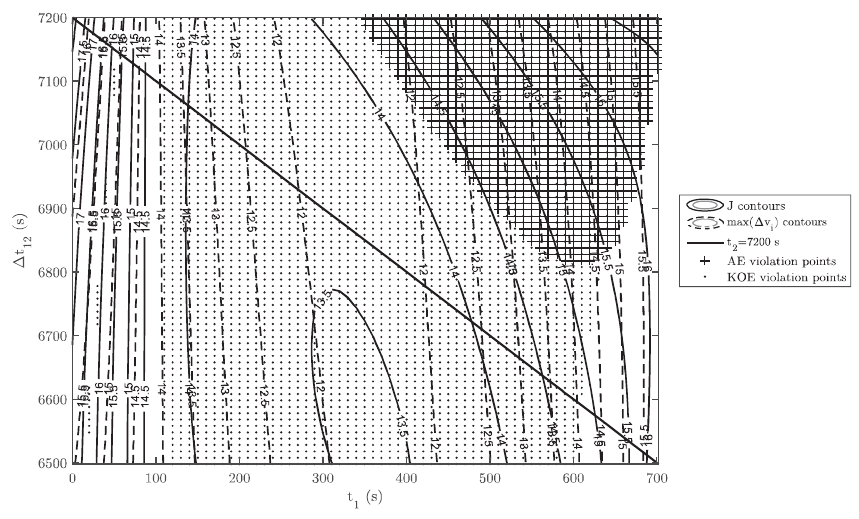}
\caption{The contours of $J$ (solid lines), the contours of $\max_i(\|\Delta\boldsymbol v_i\|)$ (dashed lines), the line of $t_2=7200\text{ s}$, and the forbidden regions according to the safety corridors.}
\label{fig:6}
\end{figure}

\begin{figure*}[!h]
\centering\includegraphics[width=0.8\linewidth]{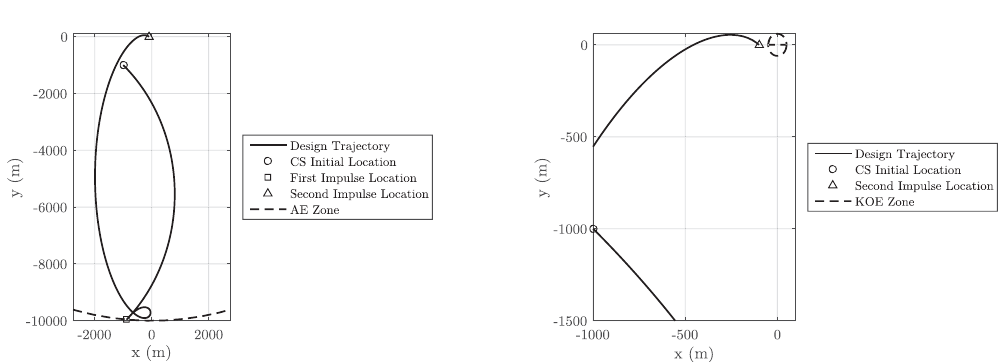}
\caption{The designed two-impulse trajectory for case study 1.}
\label{fig:7}
\end{figure*}

According to Table \ref{table:1}, from  a covariance point of view the final time should be selected as $t_2=7200\text{ s}$ that is  the maximum designated flight time. In addition, Fig. \ref{fig:4} shows that the time of the first impulse should be $t_1\leq702\text{ s}$ according to the violation of AE from $0$ to $t_1$. Since the reduction of covariance trace is the main objective in the impulse time selection, it is desired to exert the first impulse as late as possible. In this sense, according to Fig. \ref{fig:6}, $t_1^*=702\text{ s}$ and $t_2^*=7200\text{ s}$ are selected as the design points consistent with constraints and minimum covariance trace. Fig. \ref{fig:7} shows the resulting trajectory plus AE and KOE zones as well as the key locations. This trajectory requires a control effort of $J^*=16.6242\text{ m/s}$ and a maximum impulse of $\max_i⁡(\|\Delta\boldsymbol v_i\|)^*=15.9526\text{ m/s}$ at $i=1$. 

\subsubsection*{Case Study 2}

\begin{figure*}[!h]
\centering\includegraphics[width=1\linewidth]{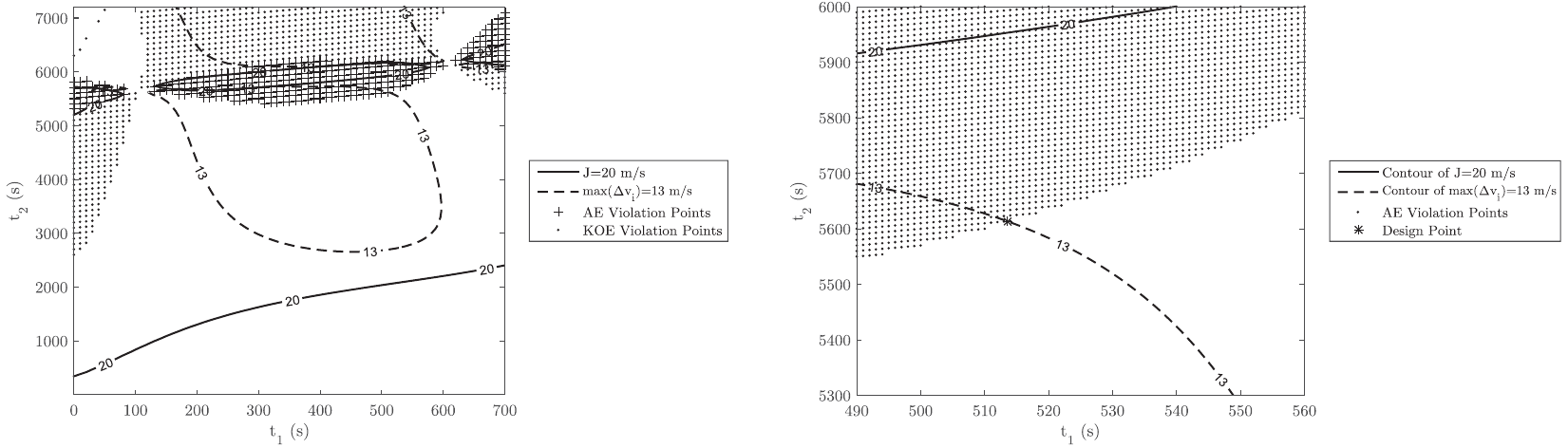}
\caption{The contours of $J=20\text{ m/s}$ and $\max_i(\|\Delta\boldsymbol v_i\|)=13\text{ m/s}$ plus the safe corridor.}
\label{fig:8}
\end{figure*}

\begin{figure*}[!h]
\centering\includegraphics[width=0.8\linewidth]{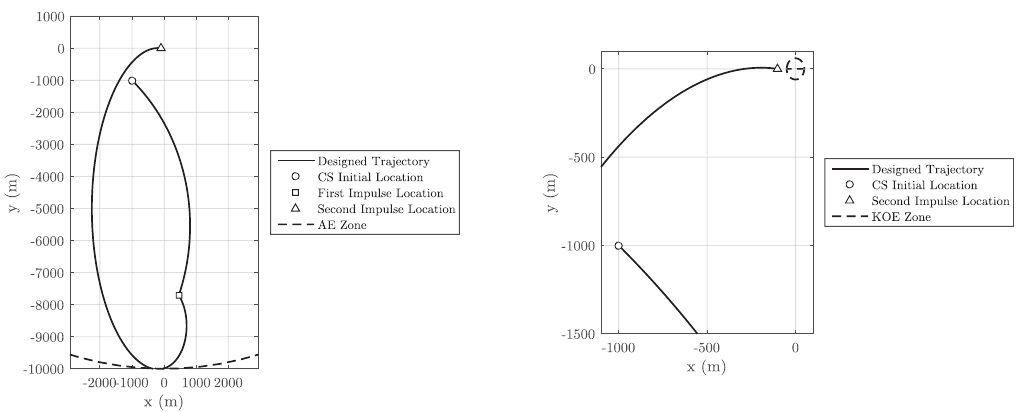}
\caption{The designed two-impulse trajectory for case study 2.}
\label{fig:9}
\end{figure*}

Fig. \ref{fig:8} shows the contours of $J=20 \text{ m/s}$, $\max⁡_i(\|\Delta\boldsymbol v_i\|)=13\text{ m/s}$, and the points where the associated trajectories violate the AE and/or the KOE after the first impulse. The maximum achievable final time occurs at the design point shown in Fig. \ref{fig:8}--right. Therefore, the design point is set for $t_1^*=513.7\text{ s}$ and $t_2^*=5613.6\text{ s}$. At this point, $J^*=14.5017\text{ m/s}$, $\max_i(\|\Delta\boldsymbol v_i\|)^*=13\text{ m/s}$, and the associated trajectory is shown in Fig. \ref{fig:9}. 

\subsection{Three-Impulse Rendezvous Design}
\label{sub:IV-B}

In a three-impulse rendezvous mission, the design parameters are $\boldsymbol r_2=\left[{r_2}_x\quad{r_2}_y\quad0\right]^T$, $t_1$, $t_2$, and $t_3$. Three of these five parameters, $t_1$, $t_2$, and $t_3$, are determined based on the characteristics of the predicted uncertainty. The maximum allowable rendezvous time duration is determined as ${t_3}_\text{max}=2\text{ h}$ according to Table \ref{table:1}. As the uncertainty covariance has a converging behavior (Fig. \ref{fig:2}), the final impulse time can be fixed at $t_3={t_3}_\text{max}=7200\text{ s}$. The first impulse is supposed to be applied at $t_1=702\text{ s}$. Parameter $t_2$ should again be selected as late as possible considering the mission constraints and the fact that some solutions for the remaining parameters, ${r_2}_x$ and ${r_2}_y$, should exist. 

\subsubsection*{Case Study 1}

The latest time to be selected for $t_2$, is $t_2=t_3$ that the three-impulse rendezvous reduces to a two-impulse rendezvous. In Subsection \ref{sub:IV-A}, for case study 1, the first impulse occurs  at $t_1=702\text{ s}$ and the final impulse occurs at the maximum allowable time. Thus, again from the covariance viewpoint in this scenario, the $N$-impulse rendezvous reduces to a two-impulse rendezvous to keep the designed key locations near the desired actual positions with highest probability. In this sense, the covariance-based three-impulse rendezvous for case study 1 occurs at $t_1^*=702\text{ s}$, $t_2^*=t_3^*=7200\text{ s}$, by selecting $\Delta\boldsymbol v_2=[0]_{3\times1}$. 

\subsubsection*{Case Study 2}

\begin{figure*}[!h]
\centering\includegraphics[width=0.9\linewidth]{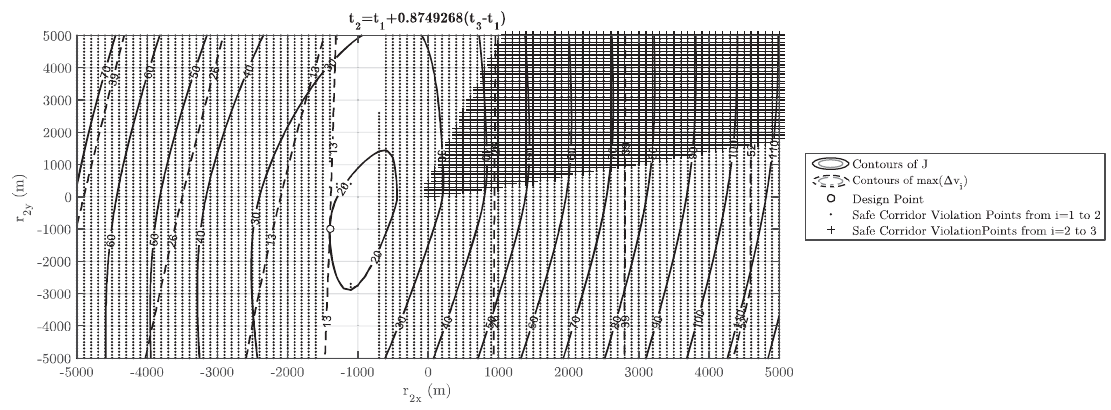}
\caption{The contours of $J$ (solid lines) and $\max_i(\|\Delta\boldsymbol v_i\|)$ (dashed lines) beside the safety corridor for case study 2.}
\label{fig:10}
\end{figure*}

\begin{figure*}[!h]
\centering\includegraphics[width=0.8\linewidth]{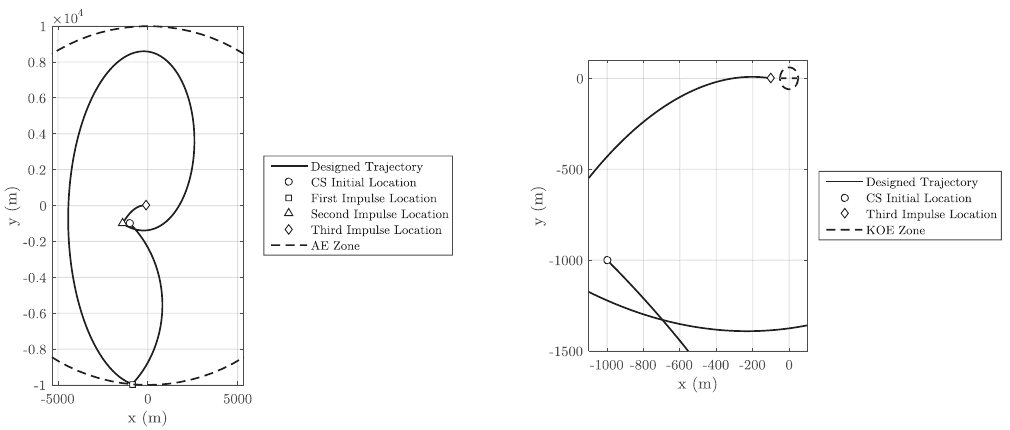}
\caption{The designed three-impulse trajectory for case study 2.}
\label{fig:11}
\end{figure*}

The second impulse time in the interval of $t_1$ to $t_3$, provides the first solution set for $\boldsymbol r_2=\left[{r_2}_x\quad{r_2}_y\quad0\right]^T$ such that the constraints get satisfied. In other words, the purpose is to find the latest time that the both loci contours of $J=20\text{ m/s}$ and $\max⁡_i(\|\Delta \boldsymbol v_i\|)=13\text{ m/s}$ lie in the safe corridor. Fig. \ref{fig:10} shows the first touch of $J=20\text{ m/s}$ and $\max⁡_i(\|\Delta\boldsymbol v_i\|)=13\text{ m/s}$ occurs at $t_2=6387 s$. So, the design point is selected to be $t_1^*=702\text{ s}$, $t_2^*=6387\text{ s}$, $t_3^*=7200\text{ s}$, and $r_2^*=-[1400.3\quad1000.0\quad0]^T\text{ m}$. The designed three-impulse trajectory has a control effort of $J^*=20.00\text{ m/s}$ and $\max_i(\|\Delta\boldsymbol v_i\|)^*=13.00\text{ m/s}$ that is shown in Fig. \ref{fig:11}.

\subsection{Discussion}
\label{sub:IV-C}

Table \ref{table:2} summarizes the results of the proposed covariance-based rendezvous trajectory design method for two simulated case studies. The following acquired remarks are notable: 
\begin{enumerate}[label=\emph{(\roman*)}]
\item In the first case study, the number of impulses has no effect on the designed trajectory. The two-impulse solution has gained the maximum impulse time and minimum uncertainties according to the safe corridor and maximum flight time. So, the $N$-impulse rendezvous problem reduces to a two impulse solution.
\item In the second case study (since the constraints are more restrictive) the two impulse solution cannot reach the minimum uncertainty at $t_1=702\text{ s}$, and $t_2=7200\text{ s}$. Thus, the final uncertainty is more than that of case study 1.
\item Applying the third impulse to the second case study, allows the spacecraft to reach its minimum uncertainty at the location of the first impulse as well as the final location.
\item In case study 2, applying the third impulse has increased the control effort as a penalty to reduce the uncertainty of the impulse locations. 
\end{enumerate}

\begin{table*}[ht]
\caption{Summarized details of the covariance-based $N$-impulse rendezvous designs.}
\begin{center}
\begin{tabular}{|c|c|c|c|c|c|c|c|c|c|}
	\hline
	Case Study & $N$ & $J\text{ (m/s)}$ & $\max_i(\|\Delta\boldsymbol v_i\|)\text{ (m/s)}$ & $\boldsymbol r_1\text{ (m)}$ & $\boldsymbol r_2\text{ (m)}$ & $\boldsymbol r_3\text{ (m)}$ & $t_1\text{ (s)}$ & $t_2\text{ (s)}$ & $t_3\text{ (s)}$ \\
	\hline
	\multirow{3}{*}{$1$} & $2$ & $16.6$ & $16.0$ & $\left(\begin{matrix} -881.0 \\ -9962.7 \end{matrix}\right)$ & $\left(\begin{matrix} -100.0 \\ 0.0 \end{matrix}\right)$ & $-$ & $702$ & $7200$ & $-$ \\
	\cdashline{2-10}
	\multicolumn{1}{|c|}{} & $3$ & $16.6$ & $16.0$ & $\left(\begin{matrix} -881.0 \\ -9962.7 \end{matrix}\right)$ & $\left(\begin{matrix} -100.0 \\ 0.0 \end{matrix}\right)$ & $-$ & $702$ & $7200$ & $-$ \\
	\hline
	\multirow{3}{*}{$2$} & $2$ & $14.5$ & $13.0$ & $\left(\begin{matrix} 460.8 \\ -7699.7 \end{matrix}\right)$ & $\left(\begin{matrix} -100.0 \\ 0.0 \end{matrix}\right)$ & $-$ & $513.7$ & $5613.6$ & $-$ \\
	\cdashline{2-10}
	\multicolumn{1}{|c|}{} & $3$ & $20.0$ & $13.0$ & $\left(\begin{matrix} -881.0 \\ -9962.7 \end{matrix}\right)$ & $\left(\begin{matrix} -1400.0 \\ -1000.0 \end{matrix}\right)$ & $\left(\begin{matrix} -100.0 \\ 0.0 \end{matrix}\right)$ & $702$ & $6387$ & $7200$ \\
	\hline
\end{tabular}
\end{center}
\label{table:2}
\end{table*}

\section{Conclusions}
\label{sec:V}

A novel methodology for space trajectory design is proposed here in which the impulse maneuvers occur at spatial points with minimum uncertainty. Covariance matrix is adopted as the measure of system uncertainty and simultaneously the solution is planned to satisfy various system constraints including maximum control effort, thruster limit, and the final time considering the safe corridors. Efficiency and viability of the proposed strategy have been investigated and verified through numerical simulations of rendezvous maneuvers exposed to various system constraints. It is demonstrated that for linear process and measurement models, selection of the impulse time instant can enhance the desired level of uncertainty in mission success, while the mission constraints can be satisfied through proper selection of the impulse positions. 

While, the covariance matrix can in general be a nonlinear function of time and/or state variables for nonlinear system and/or measurement models, future research will focus on nonlinearities associated with the covariance propagation and system complexities. Moreover, the covariance of any augmented state variables can be considered as additional mission drivers in a rendezvous mission.

\section*{Acknowledgment}

This work was supported by the Iran National Science Foundation (INSF), under grant no. 94017850, and the authors are grateful for the financial support received.

\begin{IEEEbiography}[{\includegraphics[width=1in,height=1.25in,clip,keepaspectratio]{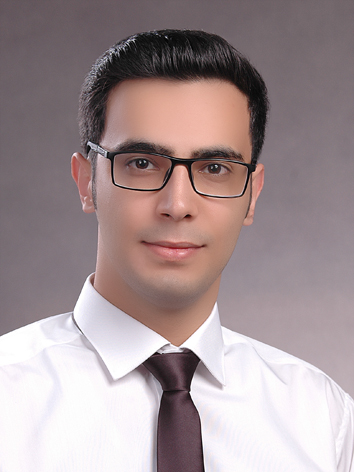}}]{Amir Shakouri}

received his MS degree in Aerospace Engineering from Sharif University of Technology (SUT), Iran, in 2017. Since 2018 he is a PhD student at SUT. His current research area of interests includes orbital maneuver, formation flying, spacecraft dynamics, and state estimation. 

\end{IEEEbiography}

\begin{IEEEbiography}[{\includegraphics[width=1in,height=1.25in,clip,keepaspectratio]{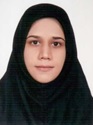}}]{Maryam Kiani}

received her PhD in Aerospace Engineering from Sharif University of Technology, Iran, in 2015. She is currently an assistant professor at Sharif university of technology, Tehran, Iran. Her main areas of interest and expertise are focused on nonlinear filtering, system identification and stochastic control with emphasis on space systems.

\end{IEEEbiography}

%\newpage

\begin{IEEEbiography}[{\includegraphics[width=1in,height=1.25in,clip,keepaspectratio]{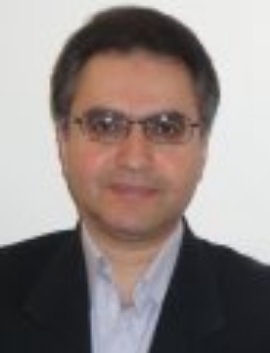}}]{Seid H. Pourtakdoust}

received the PhD degree in Aerospace Engineering from the University of Kansas, USA, in 1989. He is currently a Full Professor with the Sharif University of Technology. He has published many scientific papers. His current research interests span on aerospace flight mechanics, stochastic optimal control, attitude dynamics, and aeroelastic modeling of flapping air vehicles.

\end{IEEEbiography}

\end{document}